\documentclass[12pt]{article}
\usepackage [totalwidth=400pt, totalheight=600pt, headsep=1cm, vmarginratio=1:1]
{geometry}
\usepackage{graphicx}
\usepackage{amsmath}
\usepackage{amsfonts}
\usepackage{amssymb}
\usepackage{caption}
\usepackage[nooneline, figbotcap]{subfigure}
\usepackage[section]{placeins}
\usepackage{mathrsfs} 
\usepackage{multirow}
\usepackage{cancel}
\usepackage{yhmath}
\usepackage{pstricks-add}
\usepackage{enumerate}
\newtheorem{theorem}{Theorem}

\newtheorem{lemma}[theorem]{Lemma}

\newtheorem{problem}[theorem]{Problem}

\newtheorem{question}[theorem]{Question}

\subfiguretopcaptrue

\newcommand{\wt}[1]{
\widetriangle{#1}
}

\newcommand{\ov}[1]{
\overline{#1}
}

\newcommand{\QQ}[0]{
\mathscr{Q}
}

\newcommand{\PP}[0]{
\mathbb{P}^2
}

\newcommand{\pp}[0]{
\ell_{\infty}
}

\newcommand{\rp}[0]{
\rho_{\infty}
}

\newcommand{\EE}[0]{
\mathbb{E}^2
}

\newcommand{\RR}[0]{
\mathbb{RP}^{2}
}

\begin{document}

\title{An application of Pappus' Involution Theorem in
euclidean and non-euclidean
geometry.}

\author{Ruben Vigara \\
Centro Universitario de la Defensa - Zaragoza\\
I.U.M.A. - Universidad de Zaragoza\\}
\maketitle
\begin{abstract}
Pappus' Involution Theorem is a powerful tool for proving theorems about non-euclidean triangles and generalized triangles in Cayley-Klein models. Its power is illustrated
by proving with it some theorems about euclidean and non-euclidean polygons of different types. 
A $n$-dimensional euclidean version of these theorems is stated too.


\end{abstract}

\section{Introduction.}\label{sec:Introduction.}

In spite of being non-conformal, the use of the projective models of Cayley
\cite{Cayley} and Klein \cite{Klein} in the study of non-euclidean
planar geometries has some advantages. Any projective
theorem involving a conic could have multiple interpretations as theorems in
elliptic or hyperbolic plane. Following \cite{Richter-Gebert}, in \cite{Vigara}
those different non-euclidean
theorems emanating directly from a single projective one are called
\emph{shadows}
of the original projective theorem. In the limit case when the conic degenerates
into a single line, a non-euclidean theorem usually has a ``limit'' theorem which holds in
euclidean plane.
This property is illustrated in Section \ref{sec:generalizations}, and it has been exhaustively applied in \cite{Vigara}, where the
non-euclidean shadows of some classical projective planar theorems are explored:
the whole non-euclidean trigonometry is deduced from Menelaus' Theorem, and
Pascal's and Desargues' Theorems are used to construct some classical and
non-classical triangle centers, together with a non-euclidean version of the
Euler line and the nine-point circle of a triangle\footnote{This non-euclidean
version of the Euler line is the line denoted \emph{orthoaxis} in
\cite{Wildberger,Wildberger2}. This is the reason why it is given the name
\emph{Euler-Wildberger line} in \cite{Vigara}.}. In \cite{Vigara}, it is shown
also
the unique projective theorem hidden behind all the cosine rules of elliptic and
hyperbolic triangles and \emph{generalizes triangles} in the sense of
\cite{Buser}.

\begin{figure}
\centering
\includegraphics[width=0.5\textwidth]{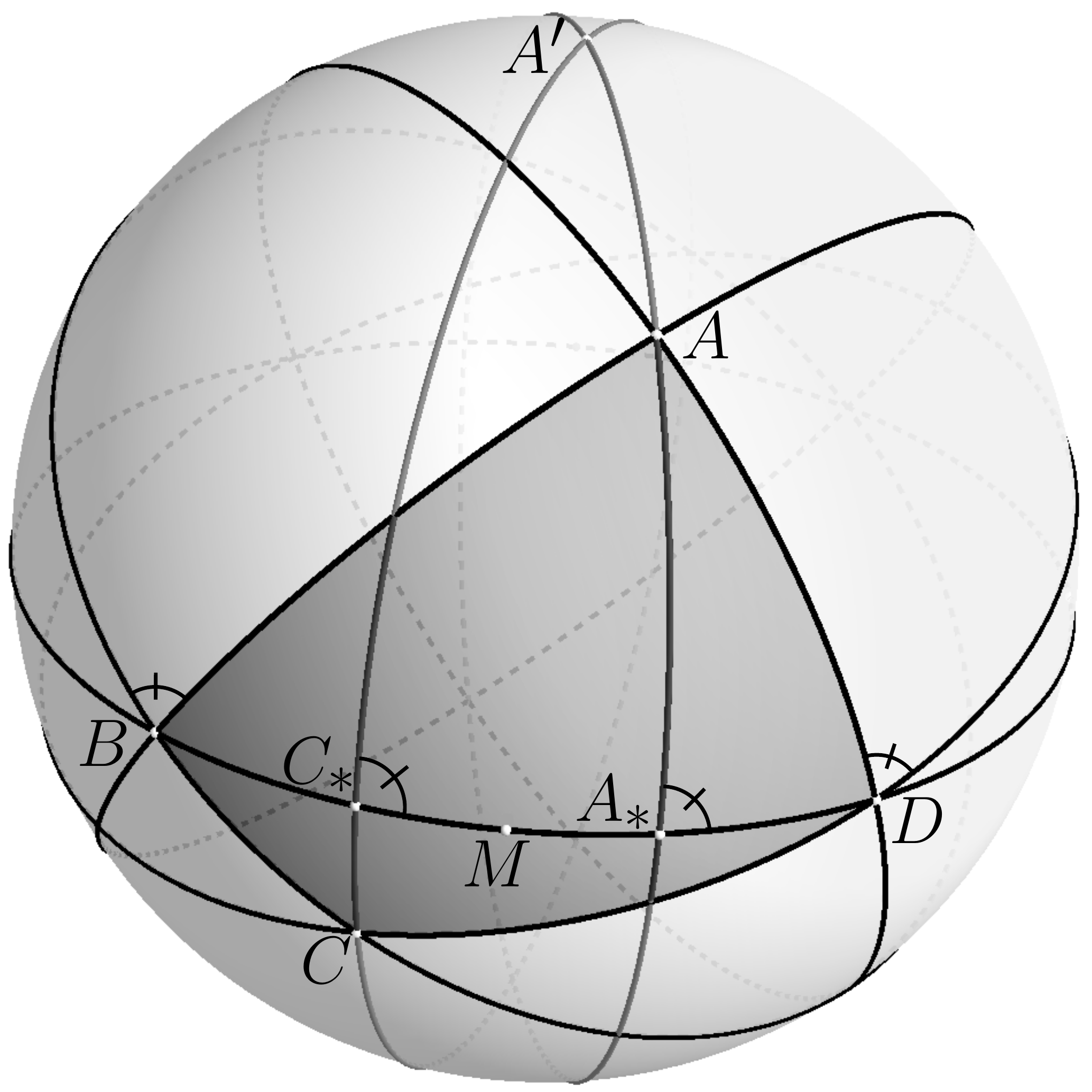}
\caption{Non-euclidean construction (spherical
view)}\label{Fig:quadrangle-spherical}
\end{figure}

Many geometric problems can be easily proven using involutions. In particular,
in many proofs and constructions of \cite{Vigara}, a
particular projective theorem has arised as an extremely
powerful tool: Pappus' Involution Theorem 
(Theorem \ref{thm:Pappus-involution} below). Here we
will exhibit its power by using it for giving a simple proof of a little
theorem about certain euclidean and non-euclidean quadrilaterals. We say that a quadrilateral in euclidean, hyperbolic 
or elliptic plane is \emph{diametral}\footnote{In euclidean plane it is a cyclic quadrilateral with a diametral diagonal.} if it has two right angles located at opposite vertices. 
\begin{theorem}\label{thm:quadrangle-absolute}
Let $R$ be a diametral quadrilateral in euclidean, hyperbolic 
or elliptic plane, with vertices $A,B,C,D$ and right angles at $B$ and $D$. 
Let $A_*,C_*$ be the
orthogonal projections of the points $A,C$ into the diagonal line $BD$,
respectively. A midpoint of the segment $\ov{BD}$ is also a midpoint of the
segment $\ov{A_*C_*}$  (see Figure
\ref{Fig:quadrangle-spherical}).
\end{theorem}
Note that in the statement of this theorem we have written ``a midpoint''
instead of ``the
midpoint''. An euclidean or hyperbolic segment is uniquely determined by its
endpoints, and it has a unique midpoint (\emph{the} midpoint). In the
elliptic case this concept is more subtle. Although we will not enter
into this discussion, depending on how we define
``segment'' and/or ``midpoint'' an elliptic segment has one or two midpoints.

We give a projective proof of Theorem \ref{thm:quadrangle-absolute} in Section
\ref{sec:proof}. The non-euclidean part of this proof relies essentially in
Pappus' Involution Theorem, and it can be reused for proving some other theorems
(shadows) about non-euclidean polygons. This will be done in Section
\ref{sec:generalizations}.
Before all of that, in Section \ref{sec:some-projective-geometry} we introduce
the basic projective tools to be used in the subsequent sections.

Finally, in Section \ref{sec:higher-dimensional} we propose a $n$-dimensional
version of Theorem \ref{thm:quadrangle-absolute},
which until now is valid only in euclidean space.

In all figures right angles are denoted with the symbol
\includegraphics[height=1em]{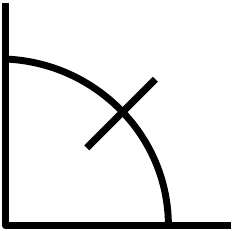}.

\section{Cayley-Klein
models}\label{sec:some-projective-geometry}

We will asume that the reader is familiar with the basic concepts
of real and complex planar projective geometry: the projective plane
and its fundamental subsets (points, lines, pencils of lines, conics),
and their projectivities. Nevertheless,
we will review some concepts and results needed for a better understanding of
Sections \ref{sec:proof} and \ref{sec:generalizations}. 
For the rigurous definitions and proofs we
refer to \cite{Cox Proj,V - Y} or \cite{Richter-Gebert}, 
for example. We assume also that the reader has
some elementary background in non-euclidean planar geometry 
(see \cite{Cox Non-euc,Santalo,Thurston}, for example).

Although we will work with real elements, we consider the real projective plane
$\mathbb{RP}^{2}$ standardly
embedded in the complex projective plane $\mathbb{CP}^{2}$. 

If $A,B$ are two different points in the projective plane we denote by $AB$ the
line joining them. If $a,b$ are two different lines, or a line and a conic, in
the projective plane, we denote by $a\cdot b$ their intersection set.

The
main theorem of projective geometry that we will use is:
\begin{theorem}[Pappus' Involution Theorem]\label{thm:Pappus-involution}
The three pairs of opposite sides of a complete quadrangle meet any line (not
through a
vertex) in three pairs of an involution.
\end{theorem}
See \cite[p. 49]{Cox Proj} for a proof. This is a partial version of 
Desargues' Involution Theorem (see \cite[p. 81]{Cox Proj}).
Using this theorem, a given complete quadrangle in the projective plane
determines a \emph{quadrangular involution} on every line not through a vertex.

Let review briefly how the projective models of euclidean, hyperbolic
and elliptic planes are constructed. We just want to show how the basic
geometric concepts needed later (perpendicular lines,
midpoint of a segment) are interpreted in projective terms, avoiding a 
full construction of these models.

\subsection{The euclidean plane}

For constructing the euclidean plane in projective terms \cite{Poncelet}, we fix a
line $\pp$ in the projective plane (the \emph{line at infinity}), and an
elliptic involution on it, i.e., a projective involution $\rp$ on $\pp$
without real fixed points (the \emph{absolute involution}). Two points
$A_\infty,B_\infty$ on $\pp$ are \emph{conjugate} if they are related by the
absolute involution: $\rp(A_\infty)=B_\infty$. The euclidean plane $\EE$ is
composed by the points of $\mathbb{RP}^{2}$ not lying in $\pp$.

For a given line
$r$ different from $\pp$, the intersection point $r\cdot \pp$ is the
\emph{point at infinity} of $r$. Two lines $r,s$ are \emph{parallel} if their
points at infinity coincide, and they are
\emph{perpendicular} if they points at infinity are conjugate.

Given two different points $A,B$ on $\EE$, the \emph{midpoint} of the segment
$\ov{AB}$ joining them is the harmonic conjugate with respect to $A,B$ of the
point at infinity of $AB$.

\begin{figure}
\centering
\includegraphics[width=\textwidth]{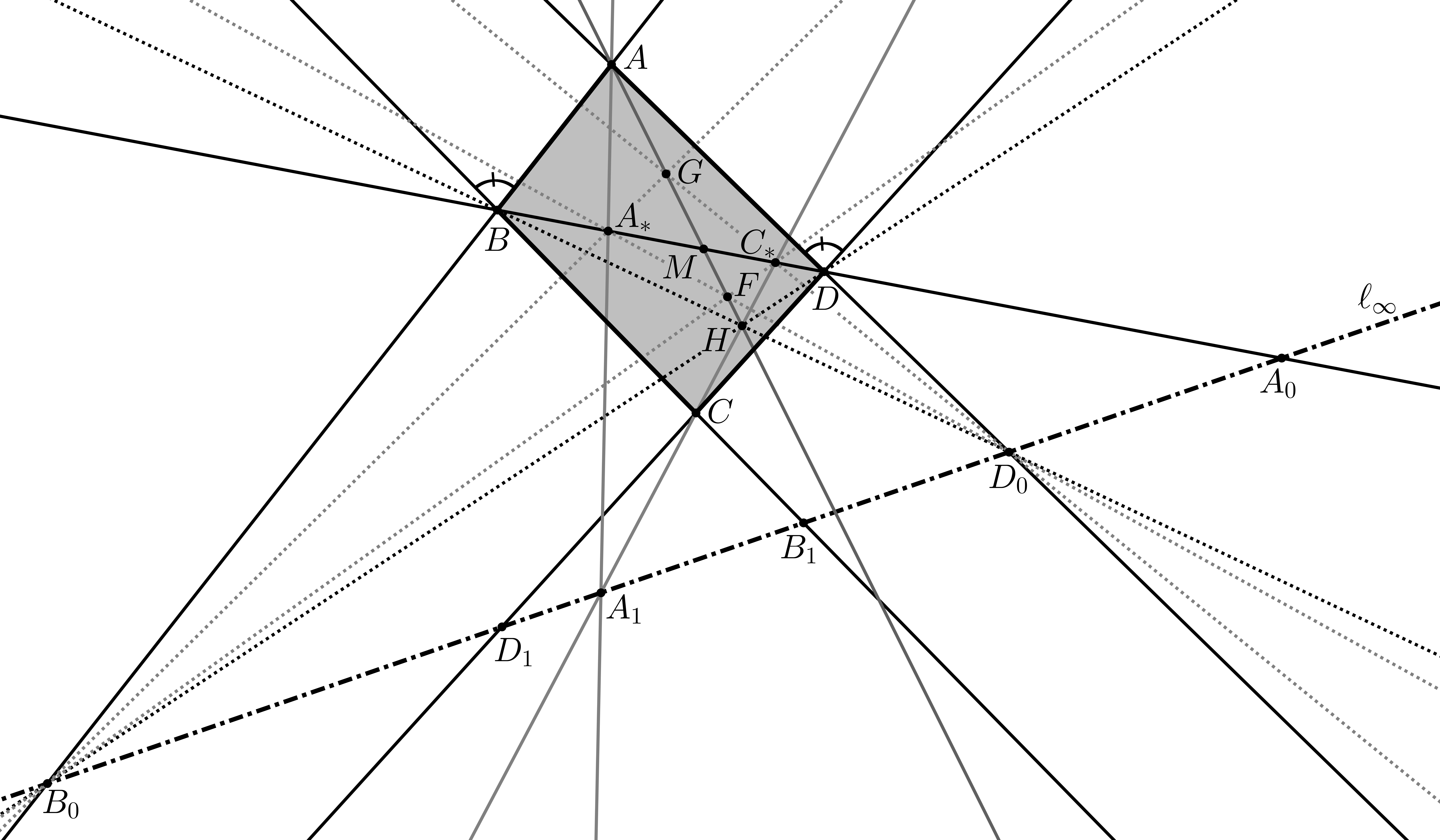}
\caption{Euclidean construction (projective
view)}\label{Fig:quadrangle-euclidean}
\end{figure}

\subsection{The hyperbolic and elliptic planes}

For constructing the non-euclidean planar models, we fix a non-degenerate conic
$\Phi_\infty$
(the \emph{absolute conic}) such that the polar of each real point with respect
to $\Phi_\infty$ is a real line. An equivalent formulation of this property is to
require, working with homogeneous coordinates, that $\Phi_\infty$ can be
expressed by an equation with real coefficients. Such a conic
can be of two kinds: a \emph{real conic}, if it has real points; or an
\emph{imaginary conic} if it has no real points (see \cite[vol. II, p. 186]{V - Y}).

When $\Phi_\infty$ is a real conic, the interior points of $\Phi_\infty$
compose the hyperbolic plane, and when $\Phi_\infty$ is an imaginary conic
the whole $\mathbb{RP}^{2}$ composes the elliptic plane.

We will use the common term \emph{the non-euclidean plane} $\PP$ either
for the hyperbolic plane (when $\Phi_\infty$ is a real conic) or for the elliptic
plane (when $\Phi_\infty$ is an imaginary conic). Geodesics in these models are given
by the intersection with $\PP$ of real projective lines. In the hyperbolic
case, we will
talk always about points or lines in a purely projective sense, even if the
referred
elements are exterior to $\Phi_\infty$.

The polarity $\rho$ with respect to $\Phi_\infty$ is a key tool in these models, 
where it plays a similar role as $\rp$ does in the euclidean case: 
two lines not tangent to $\Phi_\infty$ are \emph{perpendicular} if they are
conjugate with respect to $\Phi_\infty$, that is, if each one contains the pole of the other one 
with respect to $\Phi_\infty$. For a point $P$ and a line $p$, 
we denote by $\rho(P)$ and $\rho(p)$ the polar line of $P$ and the pole of $p$ with respect to $\Phi_\infty$, respectively. 
The polarity $\rho$ induces a natural involution 
on any line $p$ not
tangent to $\Phi_\infty$: the \emph{conjugacy involution}, which sends
each point $P\in p$ to the intersection $p\cdot \rho(P)$ of $p$ with the polar of $P$. 
The double points of the conjugacy involution on $p$ are the two points on $p\cdot \Phi_\infty$.
If the points $A,B$ not lying in $\Phi$ are conjugate with respect to $\Phi_\infty$ in the line $p$
that contains them, the polar of $A$ is the line perpendicular to $p$ through $B$ and vice versa.

Let $A,B$ be two points on $\RR$ not lying in $\Phi_\infty$ and such that the line $p$
joining them is not tangent to $\Phi_\infty$. Let $P=\rho(p)$ be the pole of $p$ with respect
to $\Phi_\infty$, and let $a,b$ be the lines joining $A,B$ with $P$, respectively.
Each of the lines $a,b$ has two (perhaps imaginary) different intersection
points with $\Phi_\infty$. Let $A_1,A_2$ and $B_1,B_2$ be the intersection points of
$a$ and $b$ with $\Phi_\infty$ respectively. The points
$$E_1=A_1B_1\cdot A_2B_2\qquad\text{and}\qquad E_2=A_1B_2\cdot A_2B_1$$
lie at the line $p$ and they are the \emph{midpoints} of the segment $\ov{AB}$.
Note that this definition is projective, and so it can be interpreted in
multiple ways (see \cite{Cox Non-euc}). For example:
\begin{itemize}
 \item If $\Phi_\infty$ is imaginary, the points $E_1,E_2$ are the two points of $p$ which are
equidistant from $A$ and $B$ in the elliptic plane $\PP$.
\item If $\Phi_\infty$ is a real conic and $A,B$ are interior
to $\Phi_\infty$, exactly one of the two points $E_1,E_2$, say $E_1$, is interior to
$\Phi_\infty$, and it is the midpoint of the hyperbolic segment $\ov{AB}$. In this
case, the other point $E_2$ is the pole of the orthogonal bisector of the
segment $\ov{AB}$.
\item If $\Phi_\infty$ is a real conic and the line $p$ is exterior to $\Phi_\infty$ its pole $P$ is interior to $\Phi_\infty$. The lines $PE_1,PE_2$ are the two bisectors of the angle between the lines $PA,PB$.
\end{itemize}

An easy characterization of midpoints is:
\begin{lemma}\label{lem:midpoints-harmonic-UV-AB}
    Let $p$ be a line not tangent to $\Phi$, let $p\cdot\Phi_\infty=\{U,V\}$,
    and take two points $A,B\in p$ different from $U,V$. 
    If $C,D$ are two points of $AB$ verifying the cross-ratio identities
	\[
	(ABCD)=(UVCD)=-1,
	\]
    then $C,D$ are the midpoints of $\ov{AB}$.
\end{lemma}

\begin{figure}
\centering
\includegraphics[width=\textwidth]{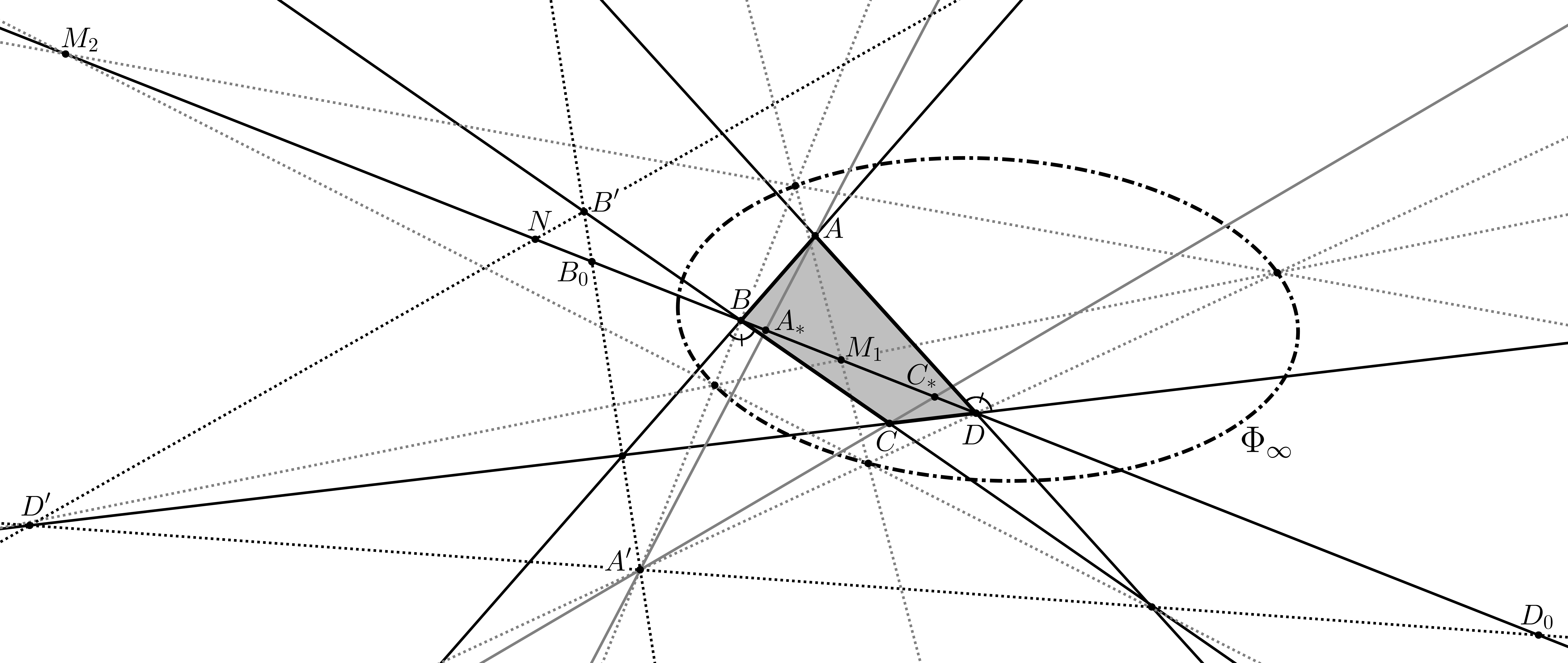}
\caption{Non-euclidean construction (hyperbolic
view)}\label{Fig:quadrangle-hyperbolic}
\end{figure}

\section{Proof of Theorem \ref{thm:quadrangle-absolute}}\label{sec:proof}

Although it is not difficult to find synthetic euclidean, elliptic or
hyperbolic proofs of Theorem \ref{thm:quadrangle-absolute}, 
we will limit ourselves to the use of projective techniques.

In both (euclidean and non-euclidean) cases, 
in the degenerate case where $A_*$ equals the point $B$ ($D$), it can be seen
that $C_*$ equals the point $D$ ($B$) and vice versa, and the statement is true.
Thus, we can assume that $A_*,C_*$ are different from $B,D$.

\paragraph{Euclidean case.}

Let $A_0,B_0,D_0$ be the points at infinity of the lines $BD, AB, AD$
respectively, and let $A_1,B_1,D_1$ be their conjugate points in $\pp$ (see
Figure \ref{Fig:quadrangle-euclidean}).
The perpendicular lines to $BD$ through $A,C$ are $AA_1,CA_1$ respectively, and
so we have
$$A_*=AA_1\cdot BD\quad \text{and}\quad C_*=CA_1\cdot BD\,.$$
Let $M$
be the midpoint of the segment $\ov{BD}$, and let $H$ be the point $BD_0\cdot
DB_0$.

If we consider the quadrangle $\QQ=\{B,C,D,H\}$, the quadrangular involution $\tau_\QQ$ that
$\QQ$ induces on $\pp$ sends $B_0,D_0$ into $B_1,D_1$ respectively, and vice
versa. This implies that $\tau_\QQ$ coincides with $\rp$ and, in particular,
that $CH$ passes through $A_1$ ($H$ is the \emph{orthocenter} of the
triangle $\wt{BCD}$). By considering the quadrangle $\{A,B_0,D_0,H\}$, it turns
out that $AH$ passes through the harmonic conjugate of $A_0$ with respect to
$B,D$, that is, that $A,H$ are collinear with $M$. In particular, this implies
that if $A_*$ and $C_*$ coincide, they coincide also with $M$.

By applying Pappus' Theorem to the hexagon $BB_0C_*A_1A_*D_0$, we have that the point
$F=B_0C_*\cdot D_0A_*$ is collinear with $A,H$. In the
same way, using the hexagon $BB_0A_*A_1C_*D_0$ it is proved that
$G=B_0A_*\cdot D_0C_*$ is collinear with $A,H$.  Taking the quadrangle
$\{F,G,B_0,D_0\}$, the point $M$ is also the harmonic conjugate of
$A_0$ with respect to $A_*,C_*$.

\begin{figure}
\centering
\includegraphics[width=0.7\textwidth]{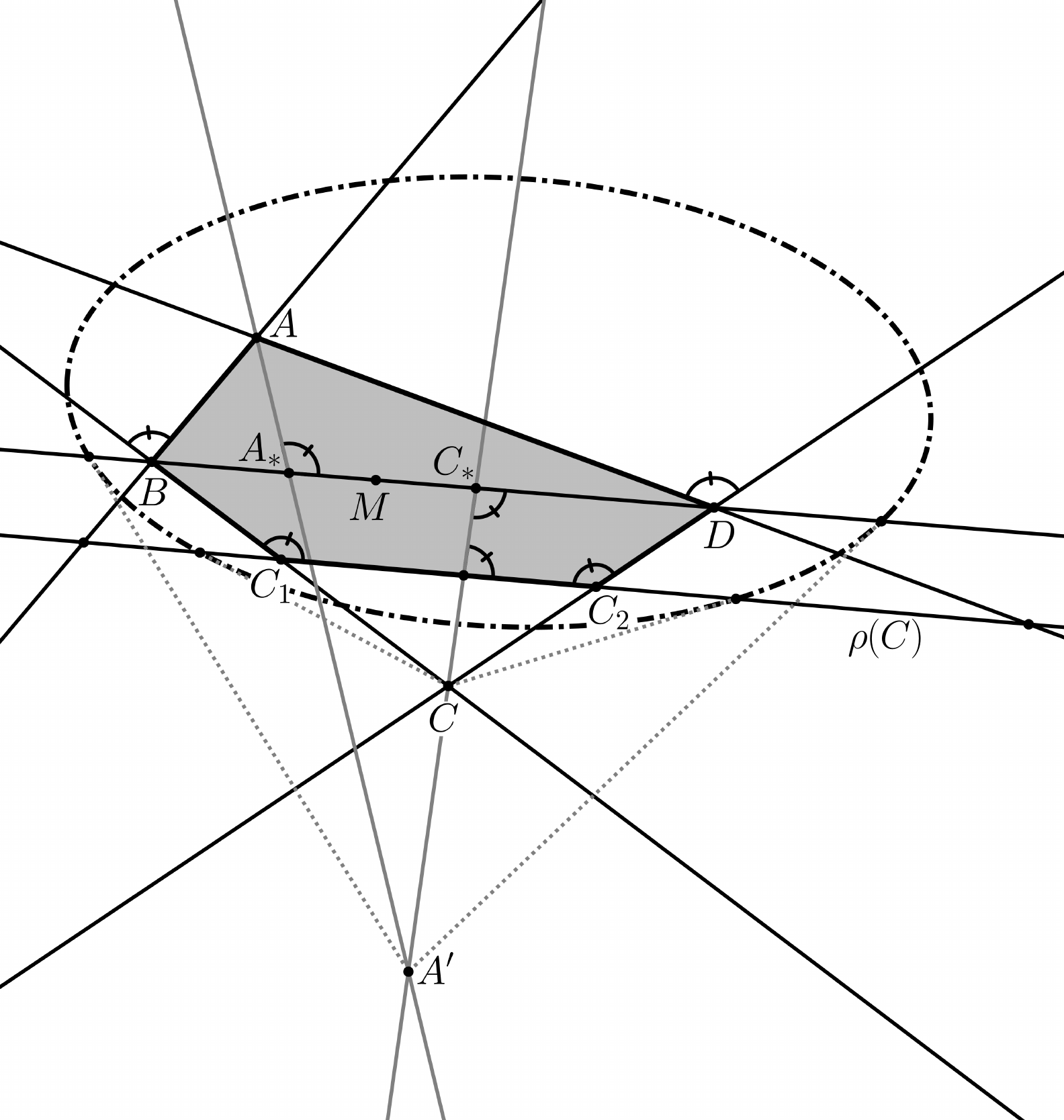}
\caption{4-right pentagon I}\label{Fig:4-right-pentagon}
\end{figure}

\paragraph{The non-euclidean case.} 
We consider the points $A,B,C,D$ such that the
lines $AB,AD$ are conjugate to $BC,DC$ respectively with respect to $\Phi_\infty$
(Figure \ref{Fig:quadrangle-hyperbolic}). 
This means that the poles $B',D'$  of the lines $AB,AD$ belong to $BC,DC$ respectively. 

Let $a$ be the line $BD$, and let $A'$ be the pole of $a$ with respect to
$\Phi_\infty$. The lines perpendicular to $a$ through $A$ and $C$ are $AA'$ and $CA'$,
respectively, and so it is $A_*=a\cdot AA'$ and $C_*=a\cdot CA'$.

Let $B_0,D_0$ be the intersection points with $a$ of the lines
$A'B',D'A'$, respectively. The triangle $\wt{A'B'D'}$ is the polar
triangle of $\wt{ABD}$, and so $B_0,D_0$ are respectively the
conjugate points of $B,D$ with respect to $\Phi_\infty$ in the line $a$. The point
$N=a\cdot B'D'$ is the pole of the line $AA'$, and so it is the conjugate point
of $A_*$ in $a$ with respect to $\Phi_\infty$.

Let $M_1,M_2$ be the midpoints of the segment $\ov{BD}$, and consider the
quadrangle $\QQ$ with vertices $C,A',B',D'$. We will make use of three involutions in $a$:
\begin{itemize}
\item the conjugacy involution $\rho_a$ induced in $a$ by the polarity with
respect to $\Phi_\infty$;
\item the quadrangular involution $\tau_\QQ$ induced in $a$ by $\QQ$; and
\item the harmonic involution $\sigma_a$ in $a$ with respect to $M_1,M_2$.
\end{itemize}

The quadrangular involution $\tau_\QQ$ sends the points $B,D,C_*$ into the points $D_0,B_0,N$ and
vice versa. This implies that the composition $\rho_a \tau_\QQ$ sends $B,D,B_0,D_0$
into $D,B,D_0,B_0$ respectively. Thus, $\sigma_a$ and $\rho_a \tau_\QQ$ agree over
at least three different points and so they coincide. As $\rho_a \tau_\QQ (C_*)=\rho_a
(N)=A_*$, the points $A_*$ and $C_*$ are harmonic conjugate with respect to $M_1$ and $M_2$,
and so by Lemma \ref{lem:midpoints-harmonic-UV-AB} the points $M_1$ and $M_2$ are the midpoints of $\ov{A_*C_*}$.

\begin{flushright}
$\blacksquare$
\end{flushright}

\begin{problem}
Find a synthetic proof of Theorem \ref{thm:quadrangle-absolute} using
the axioms of absolute geometry.
\end{problem}

\begin{figure}
\centering
\includegraphics[width=0.7\textwidth]{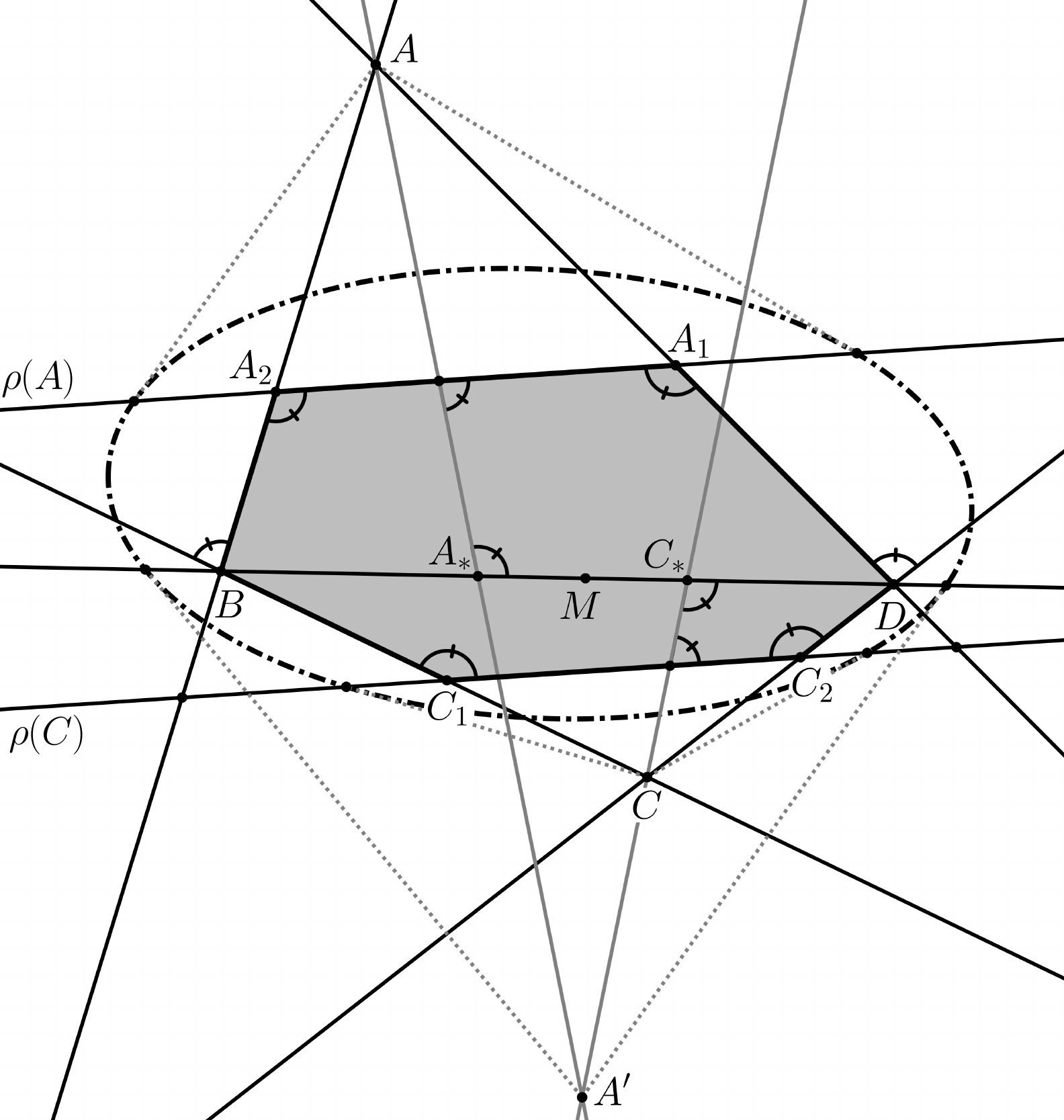}
\caption{Right-angled hexagon I}\label{Fig:right-angled-hexagon}
\end{figure}

\section{Generalizations}\label{sec:generalizations}

As we have seen, the proof of Theorem \ref{thm:quadrangle-absolute}
in the non-euclidean case is projective, and it does not depend on the
type of conic $\Phi_\infty$ that we have considered. Thus, the same proof is valid
for the hyperbolic and elliptic cases (see Figure \ref{Fig:quadrangle-spherical}). 
In the same way,
when $\Phi_\infty$ is a real conic, the same
proof does not depend on the relative position of the
vertices $A,B,C,D$ with
respect to $\Phi_\infty$.  Indeed, what we have proved in the non-euclidean
part
of the proof is the following projective theorem:
\begin{theorem}\label{thm:projective-version}
Let $\QQ=\{A,B,C,D\}$ be a complete quadrangle in the projective plane in
\emph{general position} with respect to $\Phi_\infty$ (vertices and diagonal points not
in $Phi$, sides and diagonal lines not tangent to $\Phi_\infty$) such that the lines
$AB,AD$ are conjugate to $BC,DC$ respectively with respect to $\Phi_\infty$. Let $A'$
be the pole of $BD$, and let $A_*,C_*$ be the intersection points of the lines
$AA',CA'$ respectively with $BD$. The midpoints of $\ov{A_*C_*}$ are also the
midpoints of $\ov{BD}$.
\end{theorem}

Although we were talking about quadrilaterals, Theorem
\ref{thm:projective-version} can be applied
to other hyperbolic figures that appear from the same projective
configuration.
%  We will expose some theorems for different hyperbolic polygons
% that look rather different but which are just different shadows of Theorem
% \ref{thm:projective-version}.

\begin{figure}
\centering
\includegraphics[width=0.7\textwidth]{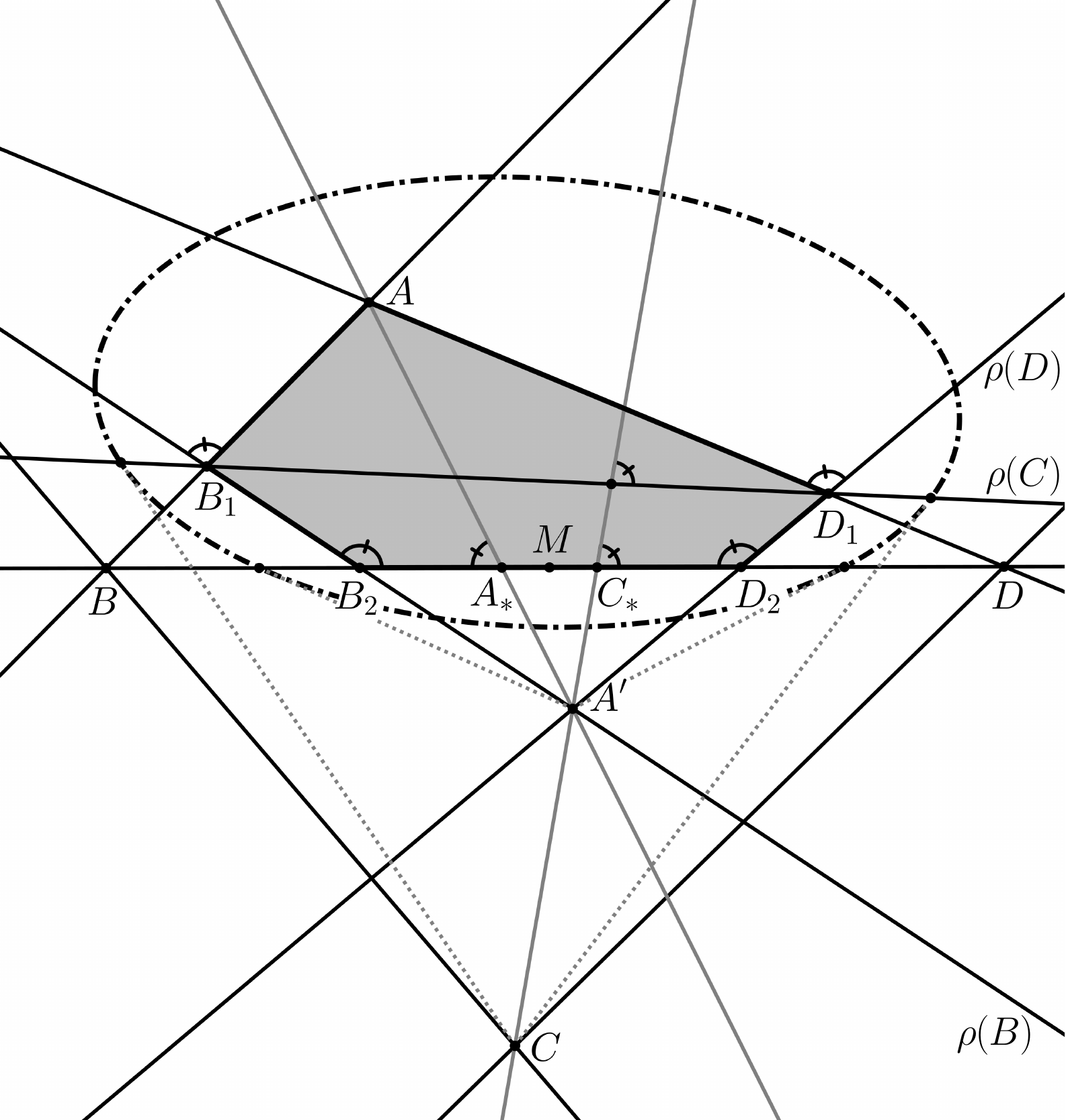}
\caption{4-right pentagon II}\label{Fig:4-right-pentagon-II}
\end{figure}

\paragraph{4-right pentagon}
If in the quadrangle $ABCD$ we assume that the
vertex $C$ lies outside the absolute conic while the rest of vertices
remain inside $\Phi_\infty$, the polar of $C$ appears into the figure as the common
perpendicular to the lines $CB$ and $CD$. The hyperbolic polygon that appears is
a
\emph{4-right pentagon}: a hyperbolic
pentagon with four right
angles (at least) at the vertices different from $A$
(Figure \ref{Fig:4-right-pentagon}).
In this case, Theorem \ref{thm:projective-version} implies:
\begin{theorem}\label{thm:pentagon-I}
In the 4-right pentagon $ABC_1C_2D$, perhaps with non-right angle at $A$, let 
$A_*$
be the orthogonal projection of $A$ into $BD$, and let $C_*$ be the intersection
of $BD$ with the common perpendicular of $BD$ and $C_1C_2$. The midpoint of
$\ov{BD}$ is also the midpoint of $\ov{A_*C_*}$.
\end{theorem}

\paragraph{Right- angled hexagon}
If in the previous figure we push also the vertex $A$ out of $\Phi_\infty$,
while $B,D$
remain interior to $\Phi_\infty$, the polar of $A$ become part of the figure as the
common perpendicular to the lines $AB$ and $AD$. The figure that appears is a
\emph{right-angled hexagon}: an hexagon in the hyperbolic plane with six right
angles as that depicted in Figure \ref{Fig:right-angled-hexagon}. With the notation of this figure, 
the traslation of Theorem
\ref{thm:projective-version} for this configuration is:
\begin{theorem}\label{thm:hexagon-I}
Let $A_1A_2BC_1C_2D$ be a right-angled hexagon. Let $A_*$ be the
intersection point with $BD$ of the common perpendicular to $BD$ and
$A_1A_2$, and let $C_*$ be the intersection point with $BD$ of the
common perpendicular to $BD$ and $C_1C_2$. The midpoint of $\ov{BD}$ is also
the midpoint of $\ov{A_*C_*}$.
\end{theorem}

\begin{figure}
\centering
\includegraphics[width=0.7\textwidth]{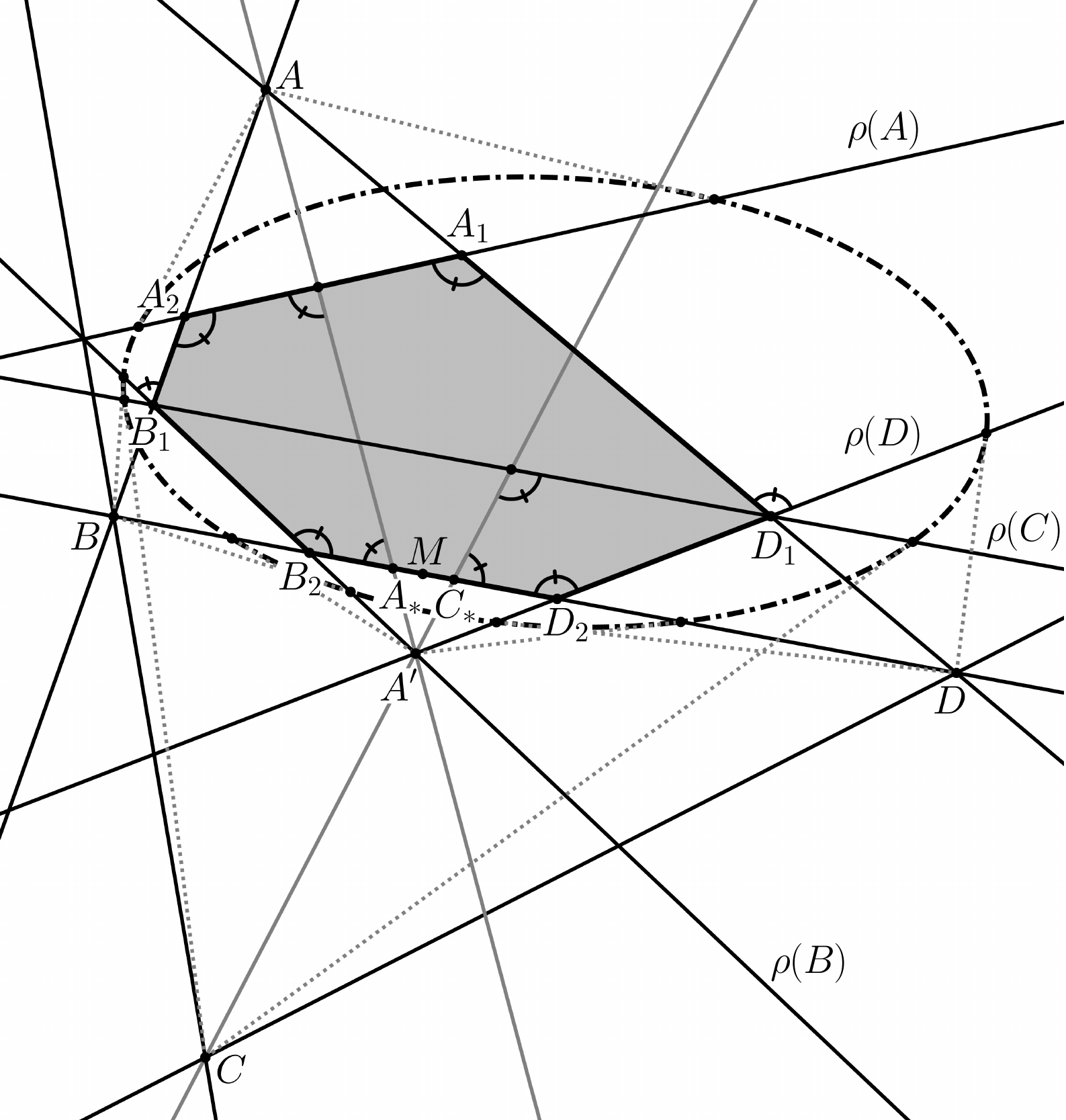}
\caption{Right-angled hexagon II}\label{Fig:right-angled-hexagon-II}
\end{figure}

\paragraph{4 right-pentagon II}
If, after pushing $C$ out of $\Phi_\infty$ for obtaining the 4-right pentagon
of Figure \ref{Fig:4-right-pentagon}, we push also
$B,D$ out of $\Phi_\infty$ but still being the line $BD$
secant to $\Phi_\infty$ ($A$ remains interior to $\Phi_\infty$), the polars of $B$ and $D$ appear
in the figure, drawing with the lines
$AB,AD$ and $BD$ another 4-right pentagon $AB_1B_2D_2D_1$ as that of
Figure \ref{Fig:4-right-pentagon-II}.
Because the pole of $AB$ is collinear with $B$ and $C$, the polars of $B$ and $C$
intersect at the point $B_1$ lying in $AB$, which is also the conjugate
point of $C_1$ in the polar of $C$ with respect to $\Phi_\infty$. In the same way, the
polars of $C$ and $D$ intersect at the point $D_1$ lying in $AD$ which is the
conjugate of $C_2$ with respect to $\Phi_\infty$. On the other hand, the polars of
$B,D$ intersect $BD$ at $B_2,D_2$ respectively, which are the conjugate points
of $B,D$ respectively in $BD$ with respect to $\Phi_\infty$. The pentagon
$AB_1B_2D_2D_1$ is a 4-right pentagon with right angles at all its vertices
with the unique possible exception of $A$. By Lemma
\ref{lem:midpoints-harmonic-UV-AB}, it can be deduced that  the midpoints of
$\ov{BD}$
are also the midpoints of $\ov{B_2D_2}$, and thus we have, 
with the notation of Figure \ref{Fig:4-right-pentagon-II}:
\begin{theorem}\label{thm:pentagon-II}
Let $AB_1B_2D_2D_1$ be a 4-right pentagon with right angles at
$B_1,B_2,D_1,D_2$. Let $A_*$ be the orthogonal projection of $A$ into
$B_2D_2$, and let $C_*$ be the intersection with $B_2D_2$ of the common
perpendicular to $B_1D_1$ and $B_2D_2$. The midpoint of $\ov{B_2D_2}$ is also
the midpoint of $\ov{A_*C_*}$.
\end{theorem}

\paragraph{Right-angled hexagon II}
If in the previous figure we push also $A$ out of the absolute conic, we obtain
a theorem similar to Theorem \ref{thm:pentagon-II} for right-angled hexagons. With the notation of Figure \ref{Fig:right-angled-hexagon-II}:
\begin{theorem}\label{thm:hexagon-II}
Let $A_1A_2B_1B_2D_2D_1$ be a right-angled hexagon. Let the common
perpendiculars to $B_2D_2$ and $A_1A_2,B_1D_1$ be denoted by $a_*,c_*$
respectively, and let $A_*,C_*$ be the intersection points of $a_*,c_*$ with
$B_2D_2$ respectively. The midpoint of $\ov{A_*C_*}$ is midpoint of
$\ov{B_2D_2}$.
\end{theorem}

\begin{figure}
\centering
\includegraphics[width=0.7\textwidth]{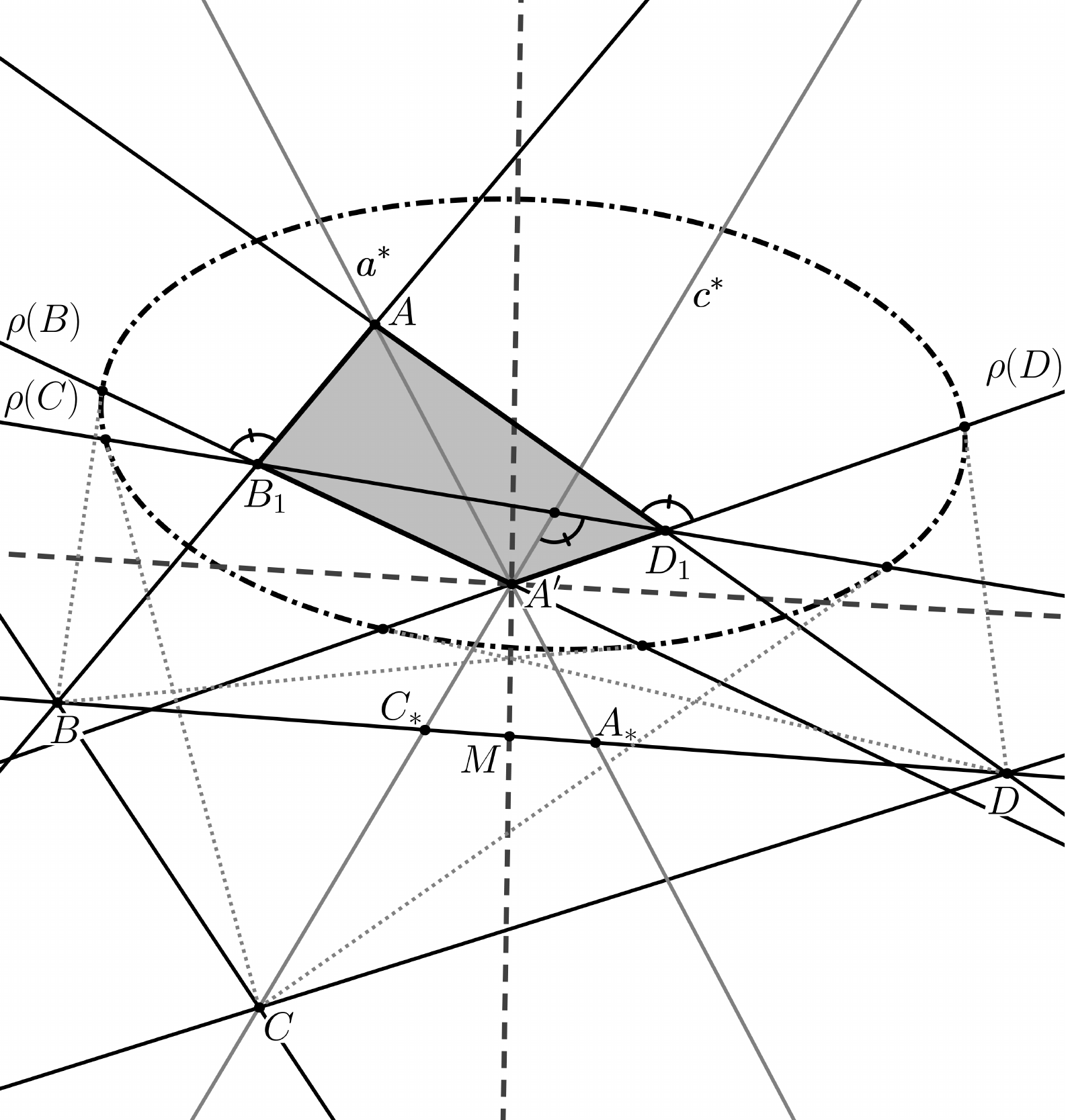}
\caption{Hyperbolic quadrangle revisited}\label{Fig:quadrangle-hyperbolic.II}
\end{figure}

\paragraph{Quadrangle II}
If in the configuration ``4 right-pentagon II'' (Figure \ref{Fig:4-right-pentagon-II})
we move $B,D$ until the line $BD$ is exterior to $\Phi_\infty$, the point $A'$
becomes interior to $\Phi_\infty$ and $AB_1A'D_1$ is a
hyperbolic diametral quadrangle with right angles at the opposite vertices $B_1,D_1$.
After a reinterpretation of the points $A_*,C_*,M_1,M_2$ for this figure, we
obtain the following theorem (see Figure \ref{Fig:quadrangle-hyperbolic.II}):
\begin{theorem}\label{thm:quadrangle-bisectors}
Let $AB_1A'D_1$ be a hyperbolic diametral quadrangle $R$ with right
angles at the opposite vertices $B_1,D_1$. Consider the lines $a_*=AA'$ and the
line $c_*$ perpendicular to $B_1D_1$. The bisectors of the angle 
$\widehat{B_1 A' D_1}$ are also the bisectors of the angle $\widehat{a_*c_*}$.
\end{theorem}
As we can expect, and as it happened with Theorem \ref{thm:quadrangle-absolute},
this theorem is also true in the euclidean and elliptic cases.
\begin{problem}\label{problem:bisectors-absolute-geometry}
Find a synthetic proof of Theorem \ref{thm:quadrangle-bisectors} using
the axioms of absolute geometry.
\end{problem}

\section{A higher-dimensional generalization}\label{sec:higher-dimensional}

A higher-dimensional generalization of Theorem \ref{thm:quadrangle-absolute} is:

\begin{theorem}\label{thm:n-simplex}
Let $\Delta$ be a simplex in euclidean $n$-dimensional space with vertices
$A_0,A_1,\ldots,A_n$. Consider the opposite face $\Delta_0$ to
$A_0$ in $\Delta$, and take the hyperplane $\pi_0$ containing $\Delta_0$. Let
$\pi_1,\pi_2,\ldots,\pi_n$ be the hyperplanes orthogonal to
$A_0A_1,A_0A_2,\ldots,A_0A_n$ through $A_1,A_2,\ldots,A_n$ respectively, and let
$C$ be the intersection point of $\pi_1,\pi_2,\ldots,\pi_n$. If $A_*,C_*$ are
the orthogonal projections of $A,C$ into $\pi_0$, the midpoint of $\ov{A_*C_*}$
is the circumcenter of $\Delta_0$.
\end{theorem}

\begin{figure}
\centering
\includegraphics[width=0.8\textwidth]{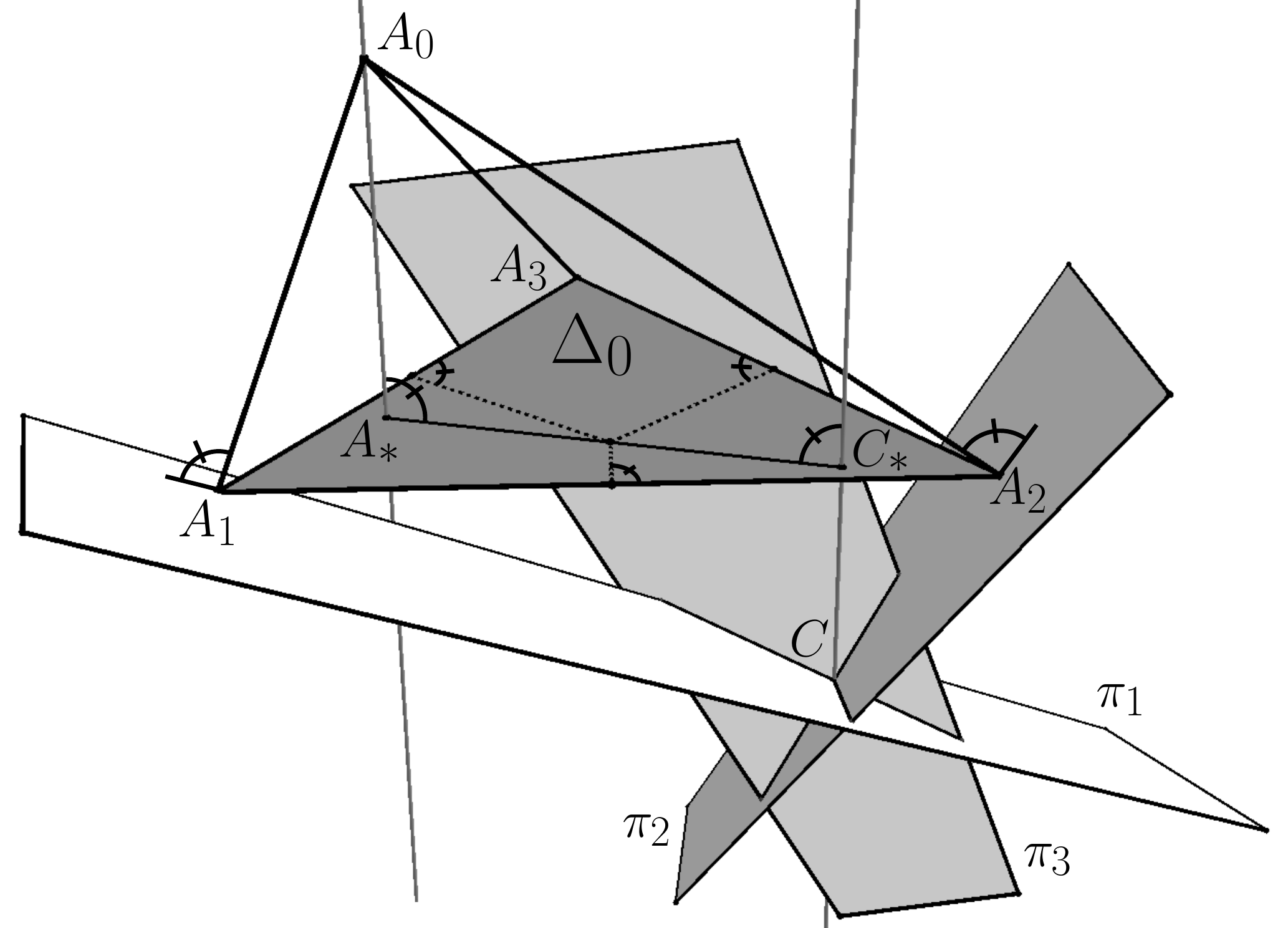}
\caption{Three-dimensional version of Theoren \ref{thm:n-simplex}}
\label{Fig:tetrahedron}
\end{figure}

We have illustrate the three-dimensional version of this theorem in Figure \ref{Fig:tetrahedron}.
Its proof (there are plenty of them) is left to the reader, it is just an
exercise on euclidean geometry. Our interest in Theorem \ref{thm:n-simplex}
relies on the fact that it is not valid in the hyperbolic and elliptic cases.

\begin{question}\label{problem:n-dimensional}
Do there exist a non-euclidean version of Theorem \ref{thm:n-simplex}?
\end{question}

According to the previous paragraph, the answer to this problem is
obviously ``no''. Nevertheless, there are many
geometric constructions that are equivalent in euclidean geometry but that are
not in the non-euclidean case. For example, in \cite{Vigara} it is shown how we can
take alternative definitions for the circumcenter and barycenter of a triangle,
different to the standard ones but equivalent to them in euclidean geometry,
in such a way that the Euler line \emph{does exist} in the hyperbolic and elliptic planes. In
Question \ref{problem:n-dimensional}, we wonder if there exists a different
formulation of Theorem \ref{thm:n-simplex} which is valid also in the
non-euclidean cases. It must be noted that in euclidean $n$-space the set of points
$A_0,A_1\ldots,A_n,C$ is \emph{diametrally cyclic}, in the sense that all these points lie in an 
$(n-1)$-dimensional sphere in $\mathbb{R}^n$ in which $A_0$ and
$C$ 
are antipodal points, while this is not the case in the non-euclidean
context.

In the same way as Theorem \ref{thm:n-simplex} is a $n$-dimensional generalization of
Theorem \ref{thm:quadrangle-absolute}, we have tried to find a $n$-dimensional 
generalization of Theorem \ref{thm:quadrangle-bisectors} without success. So the last problem that we propose is.

\begin{question}\label{problem:n-dimensional-II}
Do there exist a $n$-dimensional (euclidean or non-euclidean) version of Theorem
\ref{thm:quadrangle-bisectors}?
\end{question}

\section*{Acknowledgements}
The author wants to express his grateful thanks to Professors M. Avendano
and A.M. Oller-Marcen 
for their valuable comments and suggestions during the writing of this paper.

This research has been partially supported by the European Social Fund and 
Diputaci\'on General de Arag\'on (Grant E15 Geometr{\'\i}a).

\end{document}